\documentclass[11pt]{article}

\usepackage[margin=1in]{geometry}
\usepackage{amsmath,amssymb,amsthm,mathtools}
\usepackage{booktabs}
\usepackage{enumitem}
\setlist[itemize]{label=$\bullet$}
\usepackage{microtype}
\usepackage{xcolor}
\usepackage{graphicx}
\usepackage[colorlinks=true,linkcolor=blue!50!black,citecolor=blue!50!black,urlcolor=blue!50!black]{hyperref}

\newtheorem{theorem}{Theorem}[section]
\newtheorem{lemma}[theorem]{Lemma}
\newtheorem{proposition}[theorem]{Proposition}
\newtheorem{corollary}[theorem]{Corollary}
\theoremstyle{definition}

\theoremstyle{remark}
\newtheorem{remark}[theorem]{Remark}

\newcommand{\Zp}{\mathbb{Z}_p}
\newcommand{\Qp}{\mathbb{Q}_p}

\newcommand{\abs}[1]{\lvert #1 \rvert}
\newcommand{\norm}[1]{\lVert #1 \rVert}
\DeclareMathOperator{\vp}{v}
\newcommand{\B}{\mathrm{B}}
\newcommand{\doi}[1]{\href{https://doi.org/#1}{doi:#1}}

\title{Orbits of $p$-adic flows:\\ balls, time-one maps and minimal group actions}
\author{Kirill Terskii\\ \texttt{t02260049@gse.cs.msu.ru}}
\date{}

\begin{document}

\maketitle

\begin{abstract}
Let $f$ be a vector field on $\mathbb{Z}_p^d$ given by polynomials, or by power series converging on $\mathbb{Z}_p^d$, whose coefficients are divisible by $p^c$, where $c \ge 1$ ($c \ge 2$ when $p = 2$). Its flow exists for every $p$-adic time. I show that $t \mapsto \varphi_t(x)$ multiplies distances by exactly $\|f(x)\|_p$, in every dimension. This gives explicit flow-box coordinates, closed formulas for hitting times, return times and orbit averages, and a Newton iteration that computes them with quadratic convergence; in dimension one it recovers Mukhamedov's theorem that interpolated orbits are balls. I then show that the analytic maps $g$ congruent to the identity modulo $p^c$, coefficientwise, are exactly the time-one maps of such fields, and that the field is congruent to $g - \mathrm{id}$ modulo $p^{c+1}$. So the interpolated iterates of Bell and Poonen are flows of this kind, and I obtain a quantitative minimality criterion: if the displacements $g_i(x_0) - x_0$ of $d$ such maps span $p^k\mathbb{Z}_p^d$, the group they generate acts minimally on the ball of radius $p^{-k}$ about $x_0$. Where the argument of Cantat and Jang gives an orbit closure containing some open set, this names the ball, at every level $k$ and in every dimension; for $d = 2$ and $k = 1$ it is Jang's criterion for Markoff surfaces, which I test numerically at $7830$ points. Finally, the formal logarithm of $\mathbb{G}_m$ or of an elliptic curve is the rectifying coordinate of an explicit field, so the local divisibility $v_p([t]z) = v_p(z) + v_p(t)$ for $v_p(z) \ge 1 + c$, deep in the kernel of reduction, follows from the metric identity alone.
\end{abstract}

\section{Introduction}

Over the real numbers the flow of $dx/dt = f(x)$ is usually known only approximately, and only for a bounded time. Over the $p$-adic integers the situation is different. If $f \in p\Zp[x_1, \dots, x_d]^d$ (with $4 \mid f$ if $p = 2$), the flow $\varphi_t$ exists for every $t \in \Zp$ and is given by a power series in $t$ that converges on all of $\Zp$, and its geometry can be read off from that series. This paper is about that geometry: what the orbits look like, which maps are time-one maps of such flows, and when a group generated by such maps acts minimally on a ball.

\paragraph{The idea.}
Write the flow as $\varphi_t(x) = x + t f(x) + \sum_{n \ge 2} t^n E_n(x)$. Over $\mathbb{R}$ the higher terms are small only for small $t$, and the flow-box theorem holds on a neighbourhood of unspecified size. Over $\Zp$ there is no such trade-off. At every point $x$, each higher term is smaller than $f(x)$ by at least one $p$-adic digit, for all $t \in \Zp$ at once (Lemma~\ref{lem:key}). By the ultrametric inequality the orbit of $x$ is then, up to that one digit, the straight line $x + t f(x)$, and distances along it are exactly the distances along the line. So local questions about the flow reduce to linear algebra over $\Zp$ on the vectors $f(x)$, and, through time-one maps, questions about groups generated by maps close to the identity reduce to linear algebra on the displacements $g(x) - x$. This is why the answers below come with explicit radii, and why their hypotheses can be checked on finitely many digits of finitely many values.

\paragraph{Context.}
In arithmetic dynamics such flows usually appear as interpolations. If a map $g$ is close enough to the identity, the iterates $g^n(x)$ are the values at integers of a function analytic in $n \in \Zp$. Bell used this to prove the dynamical Mordell--Lang conjecture for automorphisms of affine varieties \cite{bell}, and Poonen gave the sharp condition on $g$ \cite{poonen}; the method is the backbone of the $p$-adic approach surveyed in \cite{bgtbook}. Cantat and Xie turned it into a $p$-adic method for actions of discrete groups \cite{cantatxie}, and Cantat and Jang used it to study orbit closures of automorphism groups of affine surfaces over $p$-adic fields \cite{cantatjang}. Jang used it to give a new proof of a result announced by Bourgain, Gamburd and Sarnak, that for the Markoff surface transitivity modulo $p$ implies minimality on the $\Zp$-points, and to extend it, for $p > 3$, to the surfaces $x^2 + y^2 + z^2 = xyz + \theta$ with $\theta \equiv 0 \pmod{p^2}$ or $\bigl(\frac{\theta - 4}{p}\bigr) = 1$ \cite[Theorem~1]{jang}. Transitivity modulo $p$ is the conjecture of Bourgain, Gamburd and Sarnak \cite{bgs}, now proved for all but finitely many $p$ \cite{chen}. In one variable, Mukhamedov showed that the interpolated orbit of a non-fixed point fills a ball \cite[Theorem~9.3]{mukhamedov}.

The iteration of a single map is also well understood: there are complete descriptions of minimal subsystems of polynomial and $1$-Lipschitz maps of $\Zp$ \cite{fanliao,anashin} and of their cycles modulo $p^n$ \cite{desjardinszieve,larin}, a local theory of analytic germs \cite{lubin,riveraletelier,jenkinsspallone,benedetto}, and, in several variables, results on periods modulo $p^n$, box dimensions of orbit closures and discrepancy \cite{mallahikaraiweiss}. Almost all of this work starts from a map, and builds a flow only when it needs one. I start from the vector field and read the geometry off its Lie series.

\paragraph{Why time is $p$-adic.}
A continuous map from $\mathbb{R}$ to $\Zp^d$ is constant, since $\mathbb{R}$ is connected and $\Zp^d$ is totally disconnected. Paths in $\Zp^d$ parametrised by real time are trivial, and the natural time domain of a $p$-adic flow is $\Zp$ itself. This changes familiar notions: time has no order, so ``the first time an event happens'' only makes sense after restricting to integer times, and the set of times at which something happens is a closed subset of $\Zp$, typically a finite union of balls.

\paragraph{Results.}
\begin{enumerate}[leftmargin=*]
\item \emph{Orbits} (Section~\ref{sec:orbits}). For a point $x$ with $f(x) \ne 0$, the map $t \mapsto \varphi_t(x)$ multiplies distances by exactly $\norm{f(x)}$, in every dimension (Theorem~\ref{thm:orbits}). In dimension one the orbit is therefore the whole ball $\B(x, \abs{f(x)})$, which is Mukhamedov's theorem for the time-one map. In higher dimension the interpolated orbit is an analytic arc in a polydisc, and Mukhamedov notes that there is then no ball-counting interpretation \cite[Section~14]{mukhamedov}; Theorem~\ref{thm:orbits} says that this arc is an isometric copy of $\Zp$, scaled by $\norm{f(x)}$. Cantat and Jang prove an analogue in normed Banach algebras, which they apply to composition operators on the Tate algebra \cite[Lemma~A.2]{cantatjang}; the pointwise statement is finer, because $\norm{f(x)}$ can be much smaller than the Gauss norm of $f$. As consequences the flow has flow-box coordinates on $\B(x, \norm{f(x)})$, first hitting times, first return times and orbit averages have closed formulas, and the integer orbit meets exactly $p^{\max(m-k,0)}$ classes modulo $p^m$, where $k = \vp(f(x))$ (Corollary~\ref{cor:count}). For time-one maps this is an exact form of the period growth $\asymp p^m$ proved for polynomial maps in \cite[Theorem~C]{mallahikaraiweiss}.
\item \emph{Time-one maps} (Section~\ref{sec:interpolated}). That an analytic map $g \equiv \mathrm{id} \pmod{p^c}$, coefficientwise, with $c$ as in \eqref{eq:hyp}, is the time-one map of a unique analytic field is the Bell--Poonen theorem \cite{bell,poonen} (in the form of \cite[Theorem~2.1]{cantatjang}; see also \cite{cantatxie}), and Cantat and Jang write the field as a logarithm \cite[(A.14)]{cantatjang}. Proposition~\ref{prop:time-one} adds the converse, so that $X \mapsto \varphi^X_1$ is a bijection from $p^c\Zp\langle x\rangle^d$ onto $\mathrm{id} + p^c\Zp\langle x\rangle^d$, and the congruence $X \equiv g - \mathrm{id} \pmod{p^{c+1}}$. The congruence is what lets the criterion below be checked on the values $g_i(x_0)$, without computing a field. As Cantat and Jang observe \cite[Remark~2.2]{cantatjang}, the congruence $g \equiv \mathrm{id}$ must hold for the coefficients and not only pointwise (see also Remark~\ref{rem:time-one-scope}).
\item \emph{Minimal actions} (Proposition~\ref{prop:minimal}). If $g_1, \dots, g_d$ are such maps and the matrix of displacements $(g_i(x_0) - x_0)/p^k$ is invertible over $\Zp$, the group they generate acts minimally on $\B(x_0, p^{-k})$, for every level $k$ and every $d$ (necessarily $k \ge c$). For $d = 2$ and $k = c = 1$ this is \cite[Proposition~4]{jang}. In the proof of \cite[Theorem~B]{cantatjang} the fields only span the tangent space over $\Qp$, and the orbit closure contains a neighbourhood of $x_0$ of unspecified size; with the $\Zp$-span I get the neighbourhood and an isometry onto it. Levels $k > 1$ cover points where no element of the group moves by exactly $p^{-1}$, the situation of \cite[Section~7.1]{jang}.
\item \emph{Markoff surfaces} (Section~\ref{sec:markoff}). I test the first-order version of the criterion numerically at $7830$ points. It holds at all but $28$ of them, and those $28$ are congruent modulo $p^2$ to the finite orbits described in \cite[Section~7.2]{jang}, in line with the conjecture stated there.
\item \emph{Examples.} Analytic fields behave like polynomial ones (Section~\ref{sec:analytic}), flow-box coordinates are computed by a Newton iteration that doubles the number of correct digits at each step (Section~\ref{sec:newton}), and the formal logarithm of $\mathbb{G}_m$ or of an elliptic curve is a flow-box coordinate of an explicit field (Section~\ref{sec:formal}).
\end{enumerate}

\paragraph{Notation.}
$\vp = \vp_p$ is the $p$-adic valuation and $\abs{x} = p^{-\vp(x)}$. On $\Zp^d$ I use $\norm{x} = \max_i \abs{x_i}$, and $\vp(x) = \min_i \vp(x_i)$. For a polynomial $g$ with coefficients in $\Zp$ (or in $\Zp^d$), $\vp(g)$ is the minimum valuation of its coefficients. $\B(x, \rho) = \{y : \norm{y - x} \le \rho\}$. For the vector field $f$ I write $c = \vp(f)$ and always assume
\begin{equation}\label{eq:hyp}
c > \tfrac{1}{p-1}, \quad\text{that is,}\quad c \ge 1 \text{ for } p \ge 3, \qquad c \ge 2 \text{ for } p = 2.
\end{equation}
A field with $p$-integral coefficients satisfies \eqref{eq:hyp} after multiplication by $p$ (by $4$ if $p = 2$), which amounts to rescaling time. The bound cannot be lowered even for linear fields: the flow $e^{\lambda t}x$ of $f(x) = \lambda x$ is defined for all $t \in \Zp$ exactly when $\vp(\lambda) > 1/(p-1)$, because the exponential series converges exactly on $\abs{z} < p^{-1/(p-1)}$. This is the rank-one case of Dwork's transfer theorem, whose proof bounds the same Taylor series $\sum_n t^n D^n/n!$ through the spectral norm of $D$ \cite{kedlaya}. Finally, $\vp(n!) = (n - s_p(n))/(p-1) \le (n-1)/(p-1)$ by Legendre's formula, where $s_p(n)$ is the sum of the base-$p$ digits of $n$.

\section{The flow and its Lie series}\label{sec:flow}

Let $D = \sum_j f_j \partial_j$ be the derivation of $\Zp[x] = \Zp[x_1, \dots, x_d]$ attached to $f$, and put
\[
E_n = \frac{D^n x}{n!} \in \Qp[x]^d, \qquad \varphi_t(x) = \sum_{n \ge 0} t^n E_n(x),
\]
where $D^n x$ is applied componentwise to the coordinate functions.

\begin{lemma}\label{lem:coeff}
$\vp(D^n g) \ge \vp(g) + nc$ for every $g \in \Zp[x]$, and $\vp(E_n) \ge nc - \vp(n!) \ge n\bigl(c - \tfrac{1}{p-1}\bigr) + \tfrac{1}{p-1}$. In particular $E_n \in \Zp[x]^d$ for all $n$, $\vp(E_n) \ge 1$ for $n \ge 1$, and $\vp(E_n) \to \infty$.
\end{lemma}

\begin{proof}
The partial derivatives preserve $\Zp[x]$ and multiplication by $f_j$ adds at least $c$ to the valuation, which gives the first claim by induction. The second follows from Legendre's formula. By \eqref{eq:hyp} the lower bound is positive for $n \ge 1$ and grows linearly.
\end{proof}

Consequently the series $\varphi_t(x)$ converges uniformly on $\Zp \times \Zp^d$ and takes values in $\Zp^d$. Write $\Phi_t = \sum_n t^n D^n/n!$ for the corresponding operator on $\Zp[x]$; by Lemma~\ref{lem:coeff} it is well defined on the Tate algebra of power series with coefficients tending to zero. Since $D$ is a derivation, $\Phi_t$ is a ring homomorphism, so $\Phi_t(g) = g \circ \varphi_t$ for polynomials $g$, and the usual identities for the exponential of a derivation hold:
\[
\Phi_{s+t} = \Phi_s \Phi_t, \qquad \frac{\partial}{\partial t} \Phi_t(g) = \Phi_t(Dg).
\]
Applied to $g = x$ they give $\varphi_{s+t} = \varphi_t \circ \varphi_s$, $\varphi_0 = \mathrm{id}$ and $\partial_t \varphi_t(x) = f(\varphi_t(x))$. Hence $\varphi$ is a flow on $\Zp^d$ with time $\Zp$, and $t \mapsto \varphi_t(x)$ is a solution given by a power series in $t$ that converges on all of $\Zp$; among such solutions it is the only one, since the coefficients of a power series solution are determined recursively by the equation.

\begin{remark}[Nilpotent fields]\label{rem:nilpotent}
For a field with $p$-integral coefficients, only the decay of $\vp(D^n x/n!)$ is needed for the Lie series to converge and define a flow, so the coefficientwise bound $c$ may be replaced by the spectral condition $\liminf_n \vp(D^n x)/n > 1/(p-1)$, as in the transfer theorem. For example $f = (p, x_1)$ has $c = 0$, but $D^2 x_1 = 0$ and $D^3 x_2 = 0$, and its flow is the polynomial map $(x_1 + pt,\ x_2 + tx_1 + pt^2/2)$. The metric results of Section~\ref{sec:orbits} use $c$ itself, and I keep \eqref{eq:hyp} throughout.
\end{remark}

\begin{remark}
Uniqueness has to be stated within analytic solutions. A $p$-adic function can have zero derivative everywhere without being constant (any locally constant function does, and there are injective examples \cite{schikhof}), so the equation alone does not determine a merely differentiable solution.
\end{remark}

\paragraph{Truncation.}
Each $\varphi_t$ with $\abs{t} \le 1$ is an isometry of $\Zp^d$: by Lemma~\ref{lem:coeff}, $\varphi_t(x) - \varphi_t(y) = (x - y) + \sum_{n \ge 1} t^n (E_n(x) - E_n(y))$ and the sum has norm at most $p^{-1}\norm{x - y}$. More importantly for computation, since $\vp(t^n E_n) \ge nc - \lfloor (n-1)/(p-1) \rfloor$ for $t \in \Zp$, let
\begin{equation}\label{eq:terms}
N(K; c) = \min \Bigl\{ N \ge 1 : nc - \bigl\lfloor \tfrac{n-1}{p-1} \bigr\rfloor \ge K \text{ for all } n \ge N \Bigr\}.
\end{equation}
The expression in $n$ is non-decreasing (its increments are $c$ or $c - 1$), so $N(K; c)$ is the first $n$ at which it reaches $K$. Then
\begin{equation}\label{eq:poly}
\varphi_t(x) \equiv \sum_{n < N(K; c)} t^n E_n(x) \pmod{p^K}, \qquad t \in \Zp, \ x \in \Zp^d.
\end{equation}
The right-hand side is a polynomial in $(t, x)$ with coefficients in $\Zp$, so $\varphi_t(x) \bmod p^K$ depends only on $t$ and $x$ modulo $p^K$. If $t \in p^m \Zp$, then $\vp(t^n E_n) \ge n(c + m) - \lfloor (n-1)/(p-1) \rfloor$ and $N(K; c + m)$ terms suffice. For example $N(19; 1) = 36$ for $p = 3$, while $N(19; 4) = 6$.

\subsection{Analytic fields}\label{sec:analytic}

Nothing above, and nothing in Section~\ref{sec:orbits}, uses that $f$ is a polynomial. Let $f$ be a vector of power series $\sum_e f_e x^e$ with $f_e \in \Zp$, $f_e \to 0$, and $\min_e \vp(f_e) = c$ satisfying \eqref{eq:hyp}; these are the elements of the Tate algebra, the series converging on $\Zp^d$. The derivation $D = f \cdot \nabla$ maps the Tate algebra to itself and raises valuations by at least $c$, so Lemmas~\ref{lem:coeff} and~\ref{lem:key} and Theorem~\ref{thm:orbits} hold verbatim.

\begin{lemma}\label{lem:analytic}
Let $\gamma = \vp((N - 1)!)$ with $N = N(K; c)$. Then $\varphi_t(x) \bmod p^K$ depends only on $f \bmod p^{K + \gamma}$.
\end{lemma}

\begin{proof}
$D^n x$ is a polynomial with integer coefficients in the coefficients of $f$ and their products with monomials; changing $f$ by a multiple of $p^{K+\gamma}$ changes $D^n x$ by a multiple of $p^{K+\gamma}$. For $n < N$, dividing by $n!$ loses at most $\gamma$ digits, and the terms $n \ge N$ vanish modulo $p^K$ by \eqref{eq:terms}.
\end{proof}

Since the coefficients of $f$ tend to zero, $f \bmod p^{K+\gamma}$ is a polynomial. So $\varphi_t(x) \bmod p^K$ for an analytic field equals the same quantity for a polynomial field, and is computed exactly by finite arithmetic.

\section{Orbits are balls}\label{sec:orbits}

The following bound is the source of all geometric statements about the flow. It sharpens Lemma~\ref{lem:coeff} at a given point.

\begin{lemma}\label{lem:key}
Let $x \in \Zp^d$ and $k = \vp(f(x))$. For $n \ge 2$,
\[
\vp\bigl(E_n(x)\bigr) \ge k + (n-1)c - \vp(n!) \ge k + 1.
\]
\end{lemma}

\begin{proof}
For $n \ge 1$ and each coordinate $i$, $D^n x_i = D(D^{n-1} x_i) = \sum_j f_j\, \partial_j (D^{n-1} x_i)$. Evaluating at $x$ and using $\vp(\partial_j D^{n-1} x_i) \ge (n-1)c$ from Lemma~\ref{lem:coeff} gives $\vp((D^n x_i)(x)) \ge k + (n-1)c$. Dividing by $n!$ gives the first inequality. For the second, $\vp(n!) \le (n-1)/(p-1) < (n-1)c$ by \eqref{eq:hyp}, and $\vp(n!)$, $(n-1)c$ are integers, so $(n-1)c - \vp(n!) \ge 1$.
\end{proof}

For $p = 2$ the hypothesis $c \ge 2$ is used here: with $c = 1$ and $n = 2$ the bound would read $k + 1 - 1 = k$, and the strict inequality fails.

\begin{theorem}\label{thm:orbits}
Let $x \in \Zp^d$ with $f(x) \ne 0$. Then for all $s, t \in \Zp$,
\[
\norm{\varphi_t(x) - \varphi_s(x)} = \norm{f(x)} \cdot \abs{t - s}.
\]
Hence $t \mapsto \varphi_t(x)$ is a homeomorphism of $\Zp$ onto the orbit of $x$, which lies in $\B(x, \norm{f(x)})$. If $d = 1$ the orbit is the whole ball $\B(x, \abs{f(x)})$.
\end{theorem}

\begin{proof}
Put $a_n = E_n(x)$, so $a_1 = f(x)$ and $\varphi_t(x) = \sum_n a_n t^n$. Then
\[
\varphi_t(x) - \varphi_s(x) = (t - s) \Bigl( a_1 + \sum_{n \ge 2} a_n h_n(s, t) \Bigr), \qquad h_n(s,t) = \sum_{i=0}^{n-1} t^i s^{n-1-i} \in \Zp.
\]
By Lemma~\ref{lem:key} each term of the sum has norm $< \norm{a_1}$, and the terms tend to zero, so the bracket has norm exactly $\norm{a_1} = \norm{f(x)}$. This proves the identity, and with it injectivity, continuity of the inverse and the inclusion in the ball (take $s = 0$).

Let $d = 1$ and $k = \vp(f(x))$, and write $F(t) = \varphi_t(x)$. For $m \ge 0$, the identity shows that $F$ induces a well-defined injective map from $\Zp/p^m\Zp$ to the set of balls of radius $p^{-k-m}$ contained in $\B(x, p^{-k})$. Both sets have $p^m$ elements, so the map is a bijection. Given $y \in \B(x, p^{-k})$, let $t_m$ be the unique class modulo $p^m$ with $F(t_m) \equiv y \pmod{p^{k+m}}$. By uniqueness the classes are compatible, so they define $t \in \Zp$, and $F(t) = y$ by continuity.
\end{proof}

\begin{remark}
In dimension one, the fact that the orbit is a ball was proved by Mukhamedov \cite[Theorem~9.3]{mukhamedov} for maps $x + p^q h(x)$ with $h \in \Zp\langle x\rangle$ and $q(p-1) > 1$, which include $\varphi_1$, through the analytic interpolation of their iterates; he writes the interpolated orbit as $a + \delta W_a(t)$, with $\delta = \varphi_1(a) - a$ and $W_a$ an analytic automorphism of the unit disc. The proof above starts from the field instead, and its metric identity holds in every dimension. The time-one map $\varphi_1$ is therefore minimal on each ball $\B(x, \abs{f(x)})$ with $f(x) \ne 0$, since $\mathbb{Z}$ is dense in $\Zp$. This agrees with the minimal decomposition theory of Fan and Liao \cite{fanliao} for polynomial maps, although $\varphi_1$ is in general not a polynomial. The setting is opposite to that of Bell and Poonen \cite{bell,poonen}, who interpolate the iterates of a given map close to the identity; Section~\ref{sec:interpolated} shows that the two settings coincide. Here the field is given, and by Mahler's theorem \cite{mahler}, $\varphi_t(x) = \sum_n \binom{t}{n} \Delta^n(x)$ with $\Delta^n(x) = \sum_{j \le n} (-1)^{n-j}\binom{n}{j} \varphi_j(x)$, so the integer iterates determine the flow at every $p$-adic time.
\end{remark}

\subsection{Flow-box coordinates}

Fix $x$ with $f(x) \ne 0$, put $k = \vp(f(x))$, choose an index $i$ with $\vp(f_i(x)) = k$, and let
\[
\B = \B(x, p^{-k}), \qquad \Sigma = \{ w \in \B : w_i = x_i \}.
\]
For $w \in \B$, $f(w) - f(x)$ has valuation at least $c + k > k$, so $\vp(f(w)) = \vp(f_i(w)) = k$ on all of $\B$, and by Theorem~\ref{thm:orbits} the orbit of every point of $\B$ stays in $\B$.

\begin{proposition}\label{prop:flowbox}
The map $\Psi(s, w) = \varphi_s(w)$ is a bijection from $\Zp \times \Sigma$ onto $\B$, and it is an isometry for the distance $\max\bigl(p^{-k}\abs{s - s'}, \norm{w - w'}\bigr)$ on $\Zp \times \Sigma$. In these coordinates $\varphi_t(\Psi(s, w)) = \Psi(s + t, w)$.
\end{proposition}

\begin{proof}
Write $\varphi_s(w) - \varphi_{s'}(w') = \bigl(\varphi_s(w) - \varphi_s(w')\bigr) + \bigl(\varphi_s(w') - \varphi_{s'}(w')\bigr) = A + C$. Since $\varphi_s$ is an isometry, $A = (w - w') + R$ with $\norm{R} \le p^{-1}\norm{w - w'}$, and $A_i = R_i$ because $w_i = w'_i$. By Theorem~\ref{thm:orbits} and $\vp(f(w')) = k$, $\norm{C} = p^{-k}\abs{s - s'}$, and the same argument applied to the $i$-th coordinate alone shows $\abs{C_i} = p^{-k}\abs{s - s'}$. If $\norm{w - w'} > p^{-k}\abs{s - s'}$, then $\norm{A + C} = \norm{w - w'}$. Otherwise $\abs{(A + C)_i} = \abs{R_i + C_i} = p^{-k}\abs{s - s'}$ (when $s \ne s'$), while $\norm{A + C} \le \max(\norm{A}, \norm{C}) = p^{-k}\abs{s - s'}$. In both cases $\norm{\Psi(s,w) - \Psi(s',w')}$ equals the stated distance, which gives isometry and injectivity.

For surjectivity, let $z \in \B$ and $G(s) = (\varphi_{-s}(z)_i - x_i)/f_i(z)$. Arguing as in the proof of Theorem~\ref{thm:orbits} for the $i$-th coordinate, $\abs{G(s) - G(s')} = \abs{s - s'}$, and $G(0) = (z_i - x_i)/f_i(z) \in \Zp$, so $G$ maps $\Zp$ to $\Zp$ isometrically. An isometry of $\Zp$ into itself is onto (count classes modulo $p^m$), so $G(s) = 0$ for some $s$, and $w = \varphi_{-s}(z) \in \Sigma$ satisfies $\Psi(s, w) = z$. The last claim is the group property of the flow.
\end{proof}

For $d = 1$ the transversal is the single point $x$, and $s(z)$ is a rectifying coordinate: $\varphi_t(z) = \varphi_{s(z) + t}(x)$.

\begin{corollary}\label{cor:times}
Let $d = 1$, $x$ with $k = \vp(f(x)) < \infty$, and $m \ge k$.
\begin{enumerate}[label=(\alph*)]
\item For $y \in \B(x, p^{-k})$ the set $\{t \in \Zp : \varphi_t(x) \in \B(y, p^{-m})\}$ is the ball $s(y) + p^{m-k}\Zp$. The first non-negative integer time at which the orbit of $x$ enters $\B(y, p^{-m})$ is the representative of $s(y)$ in $[0, p^{m-k})$. If $y \notin \B(x, p^{-k})$, the orbit never enters $\B(y, p^{-m})$.
\item The first return time of $x$ to $\B(x, p^{-m})$ under $\varphi_1$ is $p^{m-k}$.
\item The sequence $n \mapsto \varphi_n(x) \bmod p^m$ is periodic with period $p^{m-k}$ and visits each of the $p^{m-k}$ balls of radius $p^{-m}$ in $\B(x, p^{-k})$ once per period. Hence for every $F$ constant on balls of radius $p^{-m}$,
\[
\lim_{M \to \infty} \frac1M \sum_{n < M} F(\varphi_n(x)) = p^{-(m-k)} \sum_{n < p^{m-k}} F(\varphi_n(x)).
\]
\end{enumerate}
In dimension $d$, (b) and the periodicity in (c) hold with the same proof, and the orbit visits $p^{m-k}$ classes modulo $p^m$ once per period (Corollary~\ref{cor:count}). Statement (a) holds for $y$ on the orbit of $x$, with $s(y)$ defined by $\varphi_{s(y)}(x) = y$; for other $y$ one decides whether the orbit meets $\B(y, p^{-m})$ through the flow-box coordinates of Proposition~\ref{prop:flowbox}.
\end{corollary}

\begin{proof}
All three follow from $\abs{\varphi_t(x) - \varphi_{t'}(x)} = p^{-k}\abs{t - t'}$ and Theorem~\ref{thm:orbits}: $\varphi_t(x) \in \B(y, p^{-m})$ iff $\abs{t - s(y)} \le p^{-(m-k)}$.
\end{proof}

\begin{remark}[Ergodic maps and orbit closures]\label{rem:ergodic}
Let $d = 1$, $f(x_0) \ne 0$ and $k = \vp(f(x_0))$, and rescale the time-one map to $g(y) = (\varphi_1(x_0 + p^k y) - x_0)/p^k$ on $\Zp$. By Corollary~\ref{cor:times}, $g$ is an isometry of $\Zp$ that is conjugate to $y \mapsto y + 1$ through the rescaled rectifying coordinate, so it is ergodic, and each of its orbits visits the $p^m$ classes modulo $p^m$ once per period: it is a low-discrepancy sequence in the sense of \cite{mallahikaraiweiss}. Its Mahler coefficients $a_n = \Delta^n g(0)$ can be computed exactly from $g(0), \dots, g(n)$. Every $E_j$ has coefficients of valuation at least $c$, so the coefficient of $y^e$ in $g$ has valuation at least $k(e-1) + c$ for $e \ge 2$, and expanding $y^e$ in the binomial basis gives
\[
a_0 \equiv f(x_0)/p^k \not\equiv 0 \pmod p, \qquad a_1 \equiv 1 \pmod{p^c}, \qquad \vp(a_n) \ge k(n-1) + c + \vp(n!) \quad (n \ge 2).
\]
These are well within Anashin's sufficient conditions for ergodicity, as stated in \cite[Theorem~6.2]{mallahikaraiweiss}. In any dimension the orbit closures of time-one maps satisfy the dimension alternative of \cite[Theorem~B]{mallahikaraiweiss}, stated there for polynomial maps, by Corollary~\ref{cor:count}. As a pseudorandom sequence it is predictable: by \eqref{eq:poly} its terms modulo $p^K$ are a polynomial of degree less than $N(K;c)$ in $n$, so the $n$-th term costs $O(K)$ operations independently of $n$, and its index is recovered from it by Newton's method (Section~\ref{sec:newton}) in $O(\log K)$ steps. The argument uses analyticity, so it does not apply to T-functions built with OR, such as the map $x\mapsto x+(x^2\vee C)$ of \cite{klimovshamir}.
\end{remark}

\begin{corollary}[Orbit counts]\label{cor:count}
Let $x \in \Zp^d$ with $f(x) \ne 0$ and $k = \vp(f(x))$. For every $m \ge 0$ the orbit $\{\varphi_n(x) : n \ge 0\}$ meets exactly $p^{\max(m-k, 0)}$ classes modulo $p^m$. Hence its closure $\{\varphi_t(x) : t \in \Zp\}$ has box dimension one, while an equilibrium has box dimension zero. This is the dimension alternative proved for polynomial maps in \cite[Theorem~B]{mallahikaraiweiss}, here for time-one maps of analytic fields, in any dimension.
\end{corollary}

\begin{proof}
For $m \le k$ the orbit lies in $\B(x, p^{-k})$, inside one class modulo $p^m$. For $m > k$, Theorem~\ref{thm:orbits} gives $\varphi_t(x) \equiv \varphi_s(x) \pmod{p^m}$ exactly when $t \equiv s \pmod{p^{m-k}}$, and the non-negative integers meet every class modulo $p^{m-k}$. So the closure is covered by exactly $p^{\max(m-k,0)}$ balls of radius $p^{-m}$, and $\log p^{\max(m-k,0)} / \log p^m \to 1$.
\end{proof}

\subsection{Interpolated iterates and minimal group actions}\label{sec:interpolated}

Every analytic map $g \equiv \mathrm{id} \pmod{p^c}$, coefficientwise, is the time-one map of a unique analytic field. This is the Bell--Poonen theorem \cite{bell,poonen}, in the form of \cite[Theorem~2.1]{cantatjang}, and the field is the logarithm of $T_g$ applied to $\mathrm{id}$ \cite[(A.14)]{cantatjang}. I give a short proof that also yields the converse and the congruence $X \equiv g - \mathrm{id} \pmod{p^{c+1}}$, which is the part used below. For $g \in \Zp\langle x\rangle^d$ write $T_g h = h \circ g$ and $\Delta_g = T_g - \mathrm{Id}$, acting on $\Zp\langle x\rangle$ and componentwise on vectors.

\begin{proposition}[Time-one maps]\label{prop:time-one}
Let $c$ satisfy \eqref{eq:hyp}. The map $X \mapsto \varphi^X_1$ is a bijection from $p^c\Zp\langle x\rangle^d$ onto $\mathrm{id} + p^c\Zp\langle x\rangle^d$. Its inverse is
\[
X = \sum_{j \ge 1} \frac{(-1)^{j-1}}{j}\,\Delta_g^j(\mathrm{id}), \qquad X \equiv g - \mathrm{id} \pmod{p^{c+1}},
\]
and the flow of $X$ is the interpolation $\Phi(t, x) = \sum_{j \ge 0} \binom{t}{j}\Delta_g^j(\mathrm{id})(x)$ of \cite{poonen}.
\end{proposition}

\begin{proof}
If $X \in p^c\Zp\langle x\rangle^d$, then $\vp(E_n) \ge nc - \vp(n!) \ge c$ for $n \ge 1$ by Lemma~\ref{lem:coeff}, so $\varphi^X_1 \in \mathrm{id} + p^c\Zp\langle x\rangle^d$. Conversely let $g \in \mathrm{id} + p^c\Zp\langle x\rangle^d$. Since $h \circ g - h$ lies in the ideal generated by the $g_i - x_i$, $\Delta_g$ raises valuations by at least $c$, so $\vp(\Delta_g^j(\mathrm{id})) \ge jc$. Expanding $\binom{t}{j} = \sum_k s(j,k)t^k/j!$ with integer Stirling numbers $s(j,k)$ gives $\Phi(t, x) = \sum_k t^k G_k(x)$ with $G_k = \sum_{j \ge k} s(j,k)\Delta_g^j(\mathrm{id})/j!$, and $\vp(G_k) \ge \min_{j \ge k}(jc - \vp(j!)) \ge k(c - \tfrac{1}{p-1})$, so $\Phi$ converges on $\Zp \times \Zp^d$. For integers $n \ge 0$ the binomial theorem gives $\Phi(n, \cdot) = T_g^n(\mathrm{id}) = g^n$, so $\Phi(s + t, x) = \Phi(s, \Phi(t, x))$ for integers $s, t \ge 0$, and by continuity for all $s, t \in \Zp$. Both sides are power series in $s$ converging on $\Zp$, and a non-zero such series has finitely many zeros in $\Zp$ \cite{strassmann}, so their coefficients of $s$ agree: $\partial_t\Phi(t, x) = G_1(\Phi(t, x))$. Thus $t \mapsto \Phi(t, x)$ is a power series solution of $x' = G_1(x)$ with $\Phi(0, x) = x$, and by uniqueness (Section~\ref{sec:flow}) $\Phi(t, \cdot) = \varphi^{G_1}_t$, in particular $\varphi^{G_1}_1 = g$. Since $s(j,1) = (-1)^{j-1}(j-1)!$, $G_1$ is the stated series. Its $j$-th term has valuation at least $jc - \vp(j)$, which is $c$ for $j = 1$ and at least $c + 1$ for $j \ge 2$, because $\vp(j) \le \vp(j!) \le (j-1)/(p-1) < (j-1)c$ and $\vp(j)$, $(j-1)c$ are integers; this gives $G_1 \in p^c\Zp\langle x\rangle^d$ and the congruence. Finally a field is determined by its flow, $X = \partial_t\varphi^X_t|_{t=0}$, and the flow by its time-one map, since the integer iterates determine $\varphi^X_t$ by continuity; so the map is injective.
\end{proof}

Proposition~\ref{prop:time-one} is an analogue, for the group $\mathrm{id} + p^c\Zp\langle x\rangle^d$ under composition, of the correspondence between a uniform pro-$p$ group and its $\Zp$-Lie algebra given by the exponential and the logarithm \cite{ddms}. Here the Lie algebra is a module of vector fields and the exponential is the Lie series.

\begin{remark}[Which maps are flows]\label{rem:time-one-scope}
The congruence $g\equiv\mathrm{id}\pmod{p^c}$ is needed for the coefficients, not only at every point; Cantat and Jang make this point with $x+(x^q-x)^\ell$ \cite[Remark~2.2]{cantatjang}. Pointwise congruence does not even give a time-one map. The map $g(z)=z^p$ satisfies $g(z)\equiv z\pmod p$ at every $z\in\Zp$, but $g-\mathrm{id}=z^p-z$ is not divisible by $p$. It is not a time-one map: for $z\equiv w\pmod p$ the factor $\sum_iz^iw^{p-1-i}$ of $z^p-w^p$ is $\equiv pz^{p-1}\equiv0\pmod p$, so $g$ contracts each residue class, while time-one maps are isometries. Near a fixed point, Jenkins and Spallone show that $p$-adic analytic vector fields generate analytic flows and that formally equivalent maps tangent to the identity are analytically equivalent \cite{jenkinsspallone}. Crespo and Pelayo prove an existence lemma for time-dependent fields $f(x, y)$ that vanish, together with their first derivatives in $y$, along $y = v$: the solution exists for all times in a fixed ball once the initial value is close enough to $v$ \cite[Lemma~2.1]{crespopelayo}. In the language of this paper, centring at $v$ and rescaling $y - v = rz$ multiplies the terms of degree $e$ in $z$ by $r^{e-1}$, so for small $r$ the terms of degree at least two become as divisible as \eqref{eq:hyp} requires. The linear terms are not changed by rescaling, which is why they have to vanish.
\end{remark}

Thus Poonen's interpolated orbits satisfy Theorem~\ref{thm:orbits}, with $\norm{X(x)} = \norm{g(x) - x}$ by the theorem at $t = 1$, and the congruence $\Phi(t, x) \equiv x + t(g(x) - x) \pmod{p^2}$ used by Jang \cite{jang} for $c = 1$ follows from the congruence for $X$, since $\vp(E_n) \ge 2$ for $n \ge 2$ by Lemma~\ref{lem:coeff}. Flow-box coordinates for several fields at once give a criterion for minimality of the group they generate, in any dimension and on balls of any radius. For $d = 2$ and $k = 1$ it is \cite[Proposition~4]{jang}, the step through which $p$-adic flows enter the proof that residual transitivity implies minimality for Markoff surfaces.

\begin{proposition}[Minimal actions on balls]\label{prop:minimal}
Let $g_1, \dots, g_d \in \mathrm{id} + p^c\Zp\langle x\rangle^d$, $x_0 \in \Zp^d$ and $k \ge 1$, and suppose the matrix $M$ with columns $(g_i(x_0) - x_0)/p^k$ lies in $\mathrm{GL}_d(\Zp)$ (this forces $k \ge c$). Let $X_i$ be the field of $g_i$ (Proposition~\ref{prop:time-one}) and $\varphi^{(i)}$ its flow. Then
\[
F(t) = \varphi^{(1)}_{t_1} \circ \cdots \circ \varphi^{(d)}_{t_d}(x_0)
\]
is a bijection from $\Zp^d$ onto $\B = \B(x_0, p^{-k})$ with $\norm{F(t) - F(t')} = p^{-k}\norm{t - t'}$, and the group generated by $g_1, \dots, g_d$ maps $\B$ to itself and acts minimally on it.
\end{proposition}

\begin{proof}
The columns of $M$ are not divisible by $p$, so $\vp(g_i(x_0) - x_0) = k$, and by Theorem~\ref{thm:orbits} at $t = 1$, $\vp(X_i(x_0)) = k$. By Lemma~\ref{lem:key} the terms $n \ge 2$ of $g_i(x_0) - x_0 = \sum_{n \ge 1} E^{(i)}_n(x_0)$ have valuation at least $k + 1$, so $X_i(x_0) \equiv g_i(x_0) - x_0 \pmod{p^{k+1}}$. For $y \in \B$, $\vp(X_i(y) - X_i(x_0)) \ge c + k > k$, so $\vp(X_i(y)) = k$ on $\B$, and the orbits of $\varphi^{(i)}$, in particular of $g_i^{\pm1} = \varphi^{(i)}_{\pm1}$, stay in $\B$.

Let $y, y' \in \B$ and $s, s' \in \Zp$, and write
\[
\varphi^{(i)}_s(y) - \varphi^{(i)}_{s'}(y') = \bigl(\varphi^{(i)}_s(y) - \varphi^{(i)}_s(y')\bigr) + \sum_{n \ge 1}(s^n - s'^n)E^{(i)}_n(y').
\]
The first bracket is $(y - y') + R$ with $\norm{R} \le p^{-1}\norm{y - y'}$, since $\varphi^{(i)}_s$ is the identity plus a contraction (Section~\ref{sec:flow}). In the sum, the term $n = 1$ is $(s - s')X_i(x_0) + (s - s')(X_i(y') - X_i(x_0))$, and the second part has norm at most $p^{-k-1}\abs{s - s'}$. For $n \ge 2$, $\abs{s^n - s'^n} \le \abs{s - s'}$, and $\norm{E^{(i)}_n(y')} \le p^{-k-1}$ by Lemma~\ref{lem:key} at $y'$. Hence
\[
\varphi^{(i)}_s(y) - \varphi^{(i)}_{s'}(y') = (y - y') + (s - s')X_i(x_0) + e, \qquad \norm{e} \le p^{-1}\max\bigl(\norm{y - y'},\ p^{-k}\abs{s - s'}\bigr).
\]
Applying this to the $d$ factors of $F$ in turn, starting from $y = y' = x_0$, gives $F(t) - F(t') = \sum_i (t_i - t'_i)X_i(x_0) + e$ with $\norm{e} \le p^{-k-1}\norm{t - t'}$. The matrix with columns $X_i(x_0)/p^k$ is congruent to $M$ modulo $p$, hence lies in $\mathrm{GL}_d(\Zp)$, and so $\norm{F(t) - F(t')} = p^{-k}\norm{t - t'}$. Thus $F$ induces an injective map from $(\Zp/p^m)^d$ to the $p^{md}$ balls of radius $p^{-k-m}$ in $\B$, which is a bijection; as in the proof of Theorem~\ref{thm:orbits}, $F$ is onto $\B$.

For minimality let $z, z' \in \B$, and choose $t, t'$ with $F(t) = z$ and $F(t') = z'$. Then
\[
z' = \varphi^{(1)}_{t'_1} \circ \cdots \circ \varphi^{(d)}_{t'_d} \circ \varphi^{(d)}_{-t_d} \circ \cdots \circ \varphi^{(1)}_{-t_1}(z).
\]
Each $\varphi^{(i)}_s$ is an isometry and $s \mapsto \varphi^{(i)}_s(y)$ is Lipschitz uniformly in $y$, so replacing the $2d$ times by integers close to them, that is $\varphi^{(i)}_n = g_i^n$, moves the result arbitrarily little. Hence an element of the group maps $z$ arbitrarily close to $z'$.
\end{proof}

\paragraph{Relation to \cite{cantatjang,cantatxie,jang}.} The map $F$ is the map (2.10) in the proof of \cite[Theorem~B]{cantatjang}, and the density argument is theirs. There the vectors $\Theta_f(x) = X(x)$ are only required to span the tangent space over $\Qp$, which makes $F$ a local diffeomorphism. Proposition~\ref{prop:minimal} uses instead the $\Zp$-module spanned by the displacements $g(x_0) - x_0$. When this module is $p^k\Zp^d$ (by Nakayama's lemma, exactly when $d$ of the displacements form a basis of it), the minimal set is the whole ball of radius $p^{-k}$, and $F$ is an explicit isometry onto it. When such a group is a finitely generated dense subgroup of a $p$-adic analytic group $G$, acting by analytic maps congruent to the identity modulo $p$, its action extends to an analytic action of $G$ \cite{cantatxie}, and orbit closures become $G$-orbits; Proposition~\ref{prop:minimal} says which ball is a single orbit. Similarly, for character varieties of surface groups, Golsefidy and Tamam show that outside explicit exceptional sets, the closure in $\Zp$-points of an infinite orbit of the pure mapping class group contains an open set \cite[Theorem~K]{golsefidytamam}; Proposition~\ref{prop:minimal} names the ball once the displacements are known. Beyond the following numerical check I do not pursue the arithmetic of Markoff surfaces.

\subsection{A numerical check on Markoff surfaces}\label{sec:markoff}

I test the hypothesis of Proposition~\ref{prop:minimal} on the surfaces $X_\theta:\ x^2+y^2+z^2=xyz+\theta$ of \cite{bgs,jang}. Let $s_a$ be the involution that replaces $x_a$ by $x_bx_c-x_a$, where $\{a,b,c\}=\{1,2,3\}$, and let $\mathrm{rot}_a=s_b\circ s_c$. It preserves the fibres $x_a=t$ and acts on them linearly; for a point $q$ let $n_a$ be the order of this linear map modulo $p$ at $t=q_a$. Near $q$ the map $\mathrm{rot}_a^{n_a}$ is congruent to the identity on each fibre, and Proposition~\ref{prop:time-one} gives its field there.

At a smooth point $q$, two of the coordinates, centred at $q$ and divided by $p$, give an analytic chart of the residue disc of $q$ onto $\Zp^2$ by Hensel's lemma. In this chart $\mathrm{rot}_a^{n_a}$ is only affine and unipotent modulo $p$ (compare \cite[Lemma~31]{jang}; I checked this numerically at $180$ points), but its $p$-th power is congruent to the identity coefficientwise, and its field is $p$ times the field above. So if these fields and their conjugates by words in the involutions span $p\,T_q$, where $T_q=T_qX_\theta\cap\Zp^3$ is the integral tangent lattice, then Proposition~\ref{prop:minimal} with $k=1$ in the chart shows that the group acts minimally on the ball of radius $p^{-2}$ about $q$.

For $p=5,7,11$ and a range of parameters $\theta$, each replaced by a random $p$-adic integer congruent to it modulo $p^2$, I lifted every smooth point over $\mathbb F_p$ to random points $q$ and computed modulo $p^{16}$ the Smith exponents $k_1\le k_2$ of this span relative to $T_q$; the exact parameters are in the script (Section~\ref{sec:code}). Of the $7830$ lifts, $7802$ have $(k_1,k_2)=(1,1)$. The other $28$ have $(2,2)$, and each is congruent modulo $p^2$ to a point of a finite orbit, where all the fields vanish: $(0,0,\pm\sqrt\theta)$ up to permutation ($18$ lifts), a point of the Cayley cubic $X_4$ ($7$), the orbit of $(1,1,1)$ on $X_2$ ($2$) and that of $(1,\sqrt2,0)$ on $X_3$ ($1$). These are among the finite orbits of \cite[Section~7.2]{jang}, near which minimality fails. This is a computation at finitely many points, not a proof of minimality on the whole surface.

\section{Flow-box coordinates by Newton's method}\label{sec:newton}

Given $x$ and a point $z \in \B(x, p^{-k})$, I compute $(s, w) = \Psi^{-1}(z)$ from Proposition~\ref{prop:flowbox}. With $i$ as there, the iteration is
\begin{equation}\label{eq:newton}
s_0 = 0,\quad w_0 = z, \qquad \delta_j = \frac{(w_j)_i - x_i}{f_i(w_j)}, \quad s_{j+1} = s_j + \delta_j, \quad w_{j+1} = \varphi_{-\delta_j}(w_j).
\end{equation}
The division is exact: the numerator has valuation at least $k$ and the denominator exactly $k$, so $\delta_j$ is determined modulo $p^{K-k}$, which is the precision of the result.

\begin{proposition}\label{prop:newton}
Let $e_j = s - s_j$. Then $\vp(e_{j+1}) \ge 2\vp(e_j) + 1$, hence $\vp(e_j) \ge 2^j - 1$ and $\vp(\delta_j) \ge 2^j - 1$. In particular $\lceil \log_2(K - k + 1) \rceil$ steps give $s$ modulo $p^{K-k}$.
\end{proposition}

\begin{proof}
By construction $w_j = \varphi_{-s_j}(z) = \varphi_{e_j}(w)$ with $w_i = x_i$. Expanding the Lie series at $w$,
\[
(w_j)_i - x_i = e_j f_i(w) + \rho_j, \qquad \rho_j = \sum_{n \ge 2} e_j^n E_n(w)_i, \quad \vp(\rho_j) \ge k + 1 + 2\vp(e_j)
\]
by Lemma~\ref{lem:key} at $w$, where $\vp(f(w)) = k$. Also $\norm{w_j - w} = p^{-k}\abs{e_j}$ by Theorem~\ref{thm:orbits}, so $f_i(w_j) = f_i(w)(1 + \eta_j)$ with $\vp(\eta_j) \ge c + \vp(e_j) \ge 1 + \vp(e_j)$. Therefore
\[
\delta_j = \frac{e_j + \rho_j/f_i(w)}{1 + \eta_j} = e_j + O\bigl(p^{1 + 2\vp(e_j)}\bigr),
\]
and $e_{j+1} = e_j - \delta_j$ has valuation at least $2\vp(e_j) + 1$. Since $\vp(e_0) \ge 0$, induction gives $\vp(e_j) \ge 2^j - 1$, and $\delta_j = e_j - e_{j+1}$ inherits the bound.
\end{proof}

Once $s$ and $w$ are computed, Corollary~\ref{cor:times} gives first hitting times, first return times and orbit averages by closed formulas: an orbit average over balls of radius $p^{-m}$ costs $p^{m-k}$ flow evaluations instead of a long ergodic sum.

\section{The formal logarithm as a rectifying coordinate}\label{sec:formal}

The formal group of $\mathbb{G}_m$ or of an elliptic curve gives a natural example of an analytic field. Its formal logarithm is the rectifying coordinate of the flow, so Theorem~\ref{thm:orbits} contains the local divisibility $\vp([t]z) = \vp(z) + \vp(t)$ for $z \in p^{1+c}\Zp$, and Newton's method of Section~\ref{sec:newton} computes $p$-adic logarithms with no loss of digits.

Let $\mathcal F(X,Y)\in\Zp[[X,Y]]$ be a one-dimensional commutative formal group law and $\omega(z)\,dz$ its normalised invariant differential, so that $\omega\in1+z\Zp[[z]]$, $\omega(\mathcal F(X,Y))\,\partial_Y\mathcal F(X,Y)=\omega(Y)$, and
\[
\log_{\mathcal F}(z)=\int_0^z\omega=\sum_{n\ge1}\frac{\omega_{n-1}}{n}z^n,\qquad \log_{\mathcal F}(\mathcal F(X,Y))=\log_{\mathcal F}(X)+\log_{\mathcal F}(Y),
\]
and let $\exp_{\mathcal F}\in\Qp[[z]]$ be the compositional inverse of $\log_{\mathcal F}$ \cite[Sections IV.4--IV.5]{silverman}. For $\mathbb{G}_m$, $\mathcal F=X+Y+XY$, $\omega=1/(1+z)$ and $\log_{\mathcal F}(z)=\log(1+z)$. For a minimal Weierstrass equation over $\Zp$ and the parameter $z=-x/y$, the points of $E(\Qp)$ with $\vp(x)<0$, together with $O$, form the kernel of reduction $E^{(1)}(\Qp)$, and $z$ identifies it with $p\Zp$ under the group law $\mathcal F$ \cite[Propositions VII.2.1--VII.2.2]{silverman}. Write $1/\omega(z)=\sum_n w_nz^n$, with $w_n\in\Zp$ and $w_0=1$, and put
\begin{equation}\label{eq:formal-field}
f_{\mathcal F}(u)=\frac{p^c}{\omega(pu)}=\sum_{n\ge0}p^{c+n}w_nu^n,
\end{equation}
with $c$ as in \eqref{eq:hyp}. Its coefficients tend to zero and have minimum valuation $c$, so $f_{\mathcal F}$ is an analytic field in the sense of Section~\ref{sec:analytic}. For $\mathbb{G}_m$ it is the linear field $p^c(1+pu)$.

\begin{proposition}[Formal groups as flows]\label{prop:formal-flow}
For all $u,s\in\Zp$, the flow of \eqref{eq:formal-field} satisfies
\begin{equation}\label{eq:formal-translation}
p\,\varphi_s(u)=\mathcal F\bigl(pu,\exp_{\mathcal F}(p^{1+c}s)\bigr),\qquad
\log_{\mathcal F}\bigl(p\,\varphi_s(u)\bigr)=\log_{\mathcal F}(pu)+p^{1+c}s.
\end{equation}
Consequently:
\begin{enumerate}[label=(\alph*)]
\item the orbit of $0$ is $p^c\Zp$, and $s(u)=p^{-1-c}\log_{\mathcal F}(pu)$ is its rectifying coordinate;
\item $\exp_{\mathcal F}$ and $\log_{\mathcal F}$ are mutually inverse isometries of $p^{1+c}\Zp$;
\item for $z\in p^{1+c}\Zp$ and $t\in\Zp$ the multiple $[t]z=\exp_{\mathcal F}(t\log_{\mathcal F}z)$ satisfies $\vp([t]z)=\vp(z)+\vp(t)$.
\end{enumerate}
\end{proposition}

\begin{proof}
Put $e(s)=\exp_{\mathcal F}(p^{1+c}s)$. Its coefficients $p^{(1+c)n}a_n$, with $\vp(a_n)\ge-\vp(n!)$ \cite[Section IV.5]{silverman}, have valuation at least $n(1+c-\tfrac1{p-1})$, so $e$ is a power series in $s$ converging on $\Zp$, with values in $p^{1+c}\Zp$, and $\omega(e)\,e'=p^{1+c}$ since $\log_{\mathcal F}(e(s))=p^{1+c}s$. Let $\psi(s)=p^{-1}\mathcal F(pu,e(s))$. All terms of $\mathcal F(pu,e)$ have positive valuation and those of degree at least two are divisible by $p^2$, so $\psi$ is a convergent power series in $s$ with values in $\Zp$ and $\psi(0)=u$. By invariance of $\omega$,
\[
\psi'(s)=p^{-1}\,\partial_Y\mathcal F(pu,e)\,e'
=p^{-1}\frac{\omega(e)}{\omega(\mathcal F(pu,e))}\cdot\frac{p^{1+c}}{\omega(e)}
=\frac{p^c}{\omega(p\psi(s))}=f_{\mathcal F}(\psi(s)).
\]
Power series solutions are unique (Section~\ref{sec:flow}), so $\psi(s)=\varphi_s(u)$, which is the first identity of \eqref{eq:formal-translation}. Applying $\log_{\mathcal F}$, which converges on $p\Zp$, gives the second.

Since $f_{\mathcal F}(0)=p^c$, Theorem~\ref{thm:orbits} shows that the orbit of $0$ is $p^c\Zp$ and that $\abs{\varphi_t(0)-\varphi_s(0)}=p^{-c}\abs{t-s}$. By \eqref{eq:formal-translation} with $u=0$, $p\varphi_s(0)=\exp_{\mathcal F}(p^{1+c}s)$ and $\log_{\mathcal F}(p\varphi_s(0))=p^{1+c}s$. This gives (a). It also shows that $\exp_{\mathcal F}$ maps $p^{1+c}\Zp$ isometrically onto $p\cdot p^c\Zp$, with inverse $\log_{\mathcal F}$, which is (b). For (c), write $z=\exp_{\mathcal F}(w)$ with $\vp(w)=\vp(z)$; then $[t]z=\exp_{\mathcal F}(tw)$ has valuation $\vp(tw)=\vp(z)+\vp(t)$ by (b).
\end{proof}

For a point $Q$ of an elliptic curve with $\vp(z(Q))\ge1+c$, (c) is the familiar statement that $[t]Q$ lies exactly $\vp(t)$ steps deeper than $Q$ in the filtration $E^{(n)}(\Qp)=\{\vp(z)\ge n\}\cup\{O\}$ (written $E_n(\Qp)$ in \cite{silverman}). The classical proofs work in the formal group and give (b) on all of $p\Zp$ for $p$ odd \cite[Theorem IV.6.4]{silverman}; here the weaker range $p^{1+c}\Zp$ follows from the metric statement of Theorem~\ref{thm:orbits} alone. The same mechanism solves the discrete logarithm on anomalous elliptic curves: after lifting to $\Qp$ and multiplying by $p$ into the kernel of reduction, the formal logarithm turns the group law into addition \cite{smart,satoharaki}; Semaev's version of the attack works in characteristic $p$ \cite{semaev}. For $\mathbb{G}_m$ the flow is explicit, $\varphi_s(u)=\bigl((1+pu)\exp(p^{1+c}s)-1\bigr)/p$, and serves as a check.

By Proposition~\ref{prop:formal-flow}(a), Newton's method of Section~\ref{sec:newton} at $x = 0$ computes $s(u)$ modulo $p^{K-c}$ from $u \in p^c\Zp$ modulo $p^K$, hence $\log_{\mathcal F}(z)$ modulo $p^{K+1}$ from $z$ modulo $p^{K+1}$ for every $z \in p^{1+c}\Zp$, with no loss of digits.

\section{Limitations and open questions}\label{sec:limitations}

\begin{itemize}[leftmargin=*]
\item The proof of \cite[Theorem~B]{mallahikaraiweiss} appears to extend to every map in $\Zp\langle x_1, \dots, x_d\rangle^d$, which would give Corollary~\ref{cor:count} without flows; I have not checked this.
\item Proposition~\ref{prop:minimal} gives minimality on one ball. Combining it with transitivity modulo $p$ requires transitivity on the fibres modulo $p^k$, which is related to \cite[Question~32]{jang}. For the Markoff surface itself, transitivity modulo composite moduli is studied in \cite{meiripuder}, and for $p > 5$ Golsefidy and Tamam bound the number of orbits modulo $p^\ell$ uniformly in $\ell$ \cite[Proposition~N]{golsefidytamam}.
\item The computation on Markoff surfaces in Section~\ref{sec:markoff} concerns finitely many random lifts and is not a proof of minimality on the whole surface.
\item Proposition~\ref{prop:formal-flow} treats one-dimensional formal groups only.
\item Crespo and Pelayo ask for a $p$-adic version of the Eliasson--Vey linearization theorem near non-degenerate critical points of integrable systems \cite[Question~5.14]{crespopelayoww}. In the real case the flows of the quadratic parts $g_1, \dots, g_n$ of the Williamson normal form play a central role: in the elliptic case they generate a torus action, which the proofs use. Over $\Zp$ the flow of a linear field $x' = Ax$ exists for all $t \in \Zp$ as soon as the entries of $A$ are divisible by $p$ (by $4$ if $p = 2$), and this can be arranged by multiplying the $g_i$ by a power of $p$, that is, by rescaling time. Since the $g_i$ Poisson commute, their Hamiltonian flows then give an action of the compact group $\Zp^n$, which carries a Haar measure. I do not know whether these $\Zp^n$-actions can play the role of the torus action.
\end{itemize}

\section{Code availability}\label{sec:code}

The flows, flow-box coordinates, time-one maps and minimal actions are implemented in the package \texttt{galrax.padic} of the open-source library Galrax, \url{https://github.com/NotTruecCoder/galrax}, written in JAX \cite{jax} on top of Diffrax \cite{kidger}. The formal-group flows are in the file \texttt{formal.py} of that package. The numerical check on Markoff surfaces is the standalone script \texttt{anc/markoff\_smith.py} in the arXiv source of this paper. From the \texttt{anc} directory it is run as \texttt{python markoff\_smith.py} $p$ \textit{depth} \textit{lifts} $[\theta_1,\dots,\theta_r]$, where the optional list of parameters defaults to $0, \dots, p^2-1$; the three runs behind Section~\ref{sec:markoff} are listed in its docstring.

\end{document}